\documentclass[12pt,amstex]{article}

\usepackage{amssymb}
\usepackage{amsmath}
\usepackage{comment}

\def\mcM{\mathcal M}

\def\mcH{\mathcal H}

\def\s{\sigma}

\def\cav{\mbox{cav}}
\def\ve{\varepsilon}

\def\eproof{\hfill $\Box$\smallskip}

\def\beqn{\begin{equation}}
\def\eeqn{\end{equation}}

\def\endproof{\hfill $\Box$\smallskip}

\newtheorem{theorem}{Theorem}

\begin{document}

\title{The Maximal Variation of Martingales of Probabilities and Repeated Games with
Incomplete Information}

\author{Abraham Neyman\footnote{Institute Institute of
Mathematics, and Center for the Study of Rationality, The Hebrew
University of Jerusalem, Givat Ram, Jerusalem 91904, Israel. This research was supported in part
by Israel Science Foundation grants 1123/06 and 1596/10.}}
\maketitle

\begin{abstract}
The variation of a martingale $p_0^k=p_0,\ldots,p_k$ of
probabilities on a finite (or countable) set $X$ is denoted
$V(p_0^k)$ and defined by
$V(p_0^k)=E\left(\sum_{t=1}^k\|p_t-p_{t-1}\|_1\right)$. It is shown
that $V(p_0^k)\leq \sqrt{2kH(p_0)}$, where $H(p)$ is the entropy
function $H(p)=-\sum_xp(x)\log p(x)$ and $\log$ stands for the
natural logarithm. Therefore, if $d$ is the number of elements of
$X$, then $V(p_0^k)\leq \sqrt{2k\log d}$. It is shown that the order
of magnitude of the bound $\sqrt{2k\log d}$ is tight for $d\leq
2^k$: there is $C>0$ such that for every $k$ and $d\leq 2^k$ there
is a martingale $p_0^k=p_0,\ldots,p_k$ of probabilities on a set $X$
with $d$ elements, and with variation $V(p_0^k)\geq C\sqrt{2k\log
d}$. An application of the first result to game theory is that the difference
between $v_k$ and $\lim_kv_k$, where $v_k$ is the value  of the
$k$-stage repeated game with incomplete information on one side with
$d$ states,  is bounded by $\|G\|\sqrt{2k^{-1}\log
d}$ (where $\|G\|$ is the maximal absolute value of a stage payoff).
Furthermore, it is shown that the order of magnitude of this game theory bound is tight.
\end{abstract}



\section{Introduction}

Bounds on the variation of a martingale of probabilities are useful in the theory of repeated games with incomplete information. Such martingales arise as sequences of an uninformed player's posteriors $p_0^k=p_0,\ldots,p_k$ of an unknown game parameter.  The martingale's variation, $V(p_0^k):=E\sum_{t=1}^k \|p_t-p_{t-1}\|_1$, bounds from above (a positive constant times) the payoff advantage that the more informed player has over the less informed one in a two-person zero-sum $k$-stage repeated game with incomplete information on one side;  see \cite{aummas95, merzam71, merzam77, zam72}.

The maximal variation of a martingale $p_0^k$ of probabilities over a finite set depends both on the initial probability $p=p_0$, and on $k$. It is bounded by a positive constant $C(p)$ times the square root of $k$.
This inequality is used in Aumann and Maschler \cite{aummas95}\footnote{This
book is based on
 reports by Robert J. Aumann and Michael Maschler which  appeared in the sixties in
{\em Report of the U.S. Arms Control and Disarmament Agency}. See
  ``Game theoretic aspects of gradual
disarmament'' (1966, ST--80, Chapter V, pp. V1--V55), ``Repeated
games with incomplete information: a survey of recent results''
(1967, ST--116, Chapter III, pp. 287--403), and ``Repeated games with
incomplete information: the zero-sum extensive case'' (1968,
ST--143, Chapter III, pp. 37--116).} to prove that the speed of convergence of the minmax value $v_k$ of the $k$-repeated game with incomplete information on one side and perfect monitoring is $O(1/\sqrt{k})$. Zamir  \cite{zam72} proved the tightness of this bound: there is a repeated game with incomplete information on one side and perfect monitoring for which the {\em error term}, $v_k-\lim v_k$, is greater than or equal to $1/\sqrt{k}$.

Mertens and Zamir \cite{merzam71} showed that  $C(p)$ is less than or equal to $\sqrt{d-1\,}$, where $d$ is the number of elements in the support of $p$, and the error term is less than or equal to $\|G\|\sqrt{d-1\,}/\sqrt{k\,}$, where $\|G\|$ is the the maximal absolute value of a payoff in one of the possible $d$ single-stage games.

The objective of the present paper is to improve the order of magnitude of the term $\sqrt{d\,}$ in the above-mentioned bounds.
The main result of the paper is that $V(p_0^k)\leq \sqrt{2k H(p)}\leq \sqrt{2k \log d}$, where $H$ is the entropy function. This inequality implies that the error term is less than or equal to $\|G\|\sqrt{2H(p)\,}/\sqrt{k\,}$, which is less than or equal to
$\|G\|\sqrt{2\log d\,}/\sqrt{k\,}$.

We also provide tightness results for both the variation of a martingale of probabilities and the error term in repeated games with incomplete information on one side: there exists a positive constant $C$ such that for all positive integers $k$ and $d$ with $1<d\leq 2^k$ there is (1) a martingale of probabilities on a finite set with $d$ elements $p_0^k:p_0,\ldots,p_k$ with variation greater than $C\sqrt{k \log d\,}$, and (2) a repeated game with incomplete information on one side with an error term that is greater than  $C\|G\|\sqrt{\log d\,}/\sqrt{k\,}$.

\section{The results}
Let $X$ be a finite (or countable) set. For $x\in X$, the $x$-th coordinate of an element $q\in \mathbb{R}^X$ is denoted $q(x)$, and $\ell_1(X)$ is the (Banach) space of all elements $q\in \mathbb{R}^X$ with $\sum_{x\in X}|q(x)|<\infty$. Obviously, if $X$ is a finite set, then $\ell_1(X)$ equals $\mathbb{R}^X$. The $\ell_1$ norm of $q\in \ell_1(X)$ is $\|q\|_1:=\sum_{x\in X}|q(x)|$, and (thus) the $\ell_1$ distance between two elements $p,q\in \ell_1(X)$ is the sum $\|p-q\|_1=\sum_{x\in X}|p(x)-q(x)|$. A $k$-step $\ell_1(X)$-valued martingale is a stochastic process $p_0^k=p_0,\ldots,p_k$ where $p_t$, $0\leq t\leq k$, takes values in $\ell_1(X)$ and $E(p_t\mid p_0,\ldots,p_{t-1})=p_{t-1}$. Let $\Delta(X)$ denote all probabilities on $X$, i.e., all elements $p\in \mathbb{R}^X_+$ with $\sum_{x\in X}p(x)=1$, and for $p\in \Delta(X)$ and a
positive integer $k$ we denote by $\mcM_k(X,p)$ the set of all
martingales $p_0^k$ with $p_t\in \Delta(X)$ and $p_0=p$.

The variation
of the martingale $p_0^k$ is denoted $V(p_0^k)$ and is defined by
$V(p_0^k)=E\left(\sum_{t=1}^k\|p_t-p_{t-1}\|_1\right)$.
Set
\beqn V(k,p):=\sup\{V(p_0^k): p_0^k\in \mcM_k(X,p)\} \eeqn
and
\beqn V(k,d):=\sup\{V(k,p): p\in \Delta(X)
\mbox{ and } |X|=d\}. \eeqn

A trivial inequality is $V(k,p)\leq 2k$. A classical bound (that is used in the theory of repeated games with
incomplete information; see \cite{aummas95, merzam71}) of $V(k,d)$ is
\[V(k,d)\leq \sqrt{k(d-1)}.\]
This classical bound improves the trivial bound only for $d\leq 4k$.
Our objective is to derive a meaningful bound that (1) is applicable
also to $d>4k$, and (2) such that its order of magnitude is the best
possible one for large $d$. We have

\begin{theorem}\label{thm:martingalevariationbound}
\[V(k,p)\leq \sqrt{2k H(p)}\]
and thus
\[V(k,d)\leq \sqrt{2k\log d}\,,\]
where $H(p)=-\sum_x p(x)\log p(x)$ is the entropy function and $\log$ stands for the natural
logarithm.
\end{theorem}

As $V(k,p)\leq 2k$, the results of Theorem \ref{thm:martingalevariationbound} are of interest for
$H(p)\leq 2k$ and for $d\leq e^{2k}$. For large values of $d\leq
e^{2k}$, the bound $\sqrt{2k\log d}$ is a significant improvement
over the classical bound $\sqrt{k(d-1)}$. Moreover, as there are
probabilities $p$ over a countable set $X$ with finite entropy, the
bound $\sqrt{2kH(p)}$ is applicable independently of the size of the
set $X$.

One may wonder if the order of magnitude of each of the bounds,
$\sqrt{2k\log d}$ and $\sqrt{2k H(p)}$, are the best possible. For
$X=\{0,1\}$ and $p(\alpha)=(\alpha,1-\alpha)\in \Delta(X)$ we have
$V(k,p(\alpha))\leq \sqrt{k\alpha(1-\alpha)}$. As
$\alpha(1-\alpha)=o(H(p(\alpha)))$ as $\alpha\to 0+$, the order of
magnitude of the bound $\sqrt{k H(p)}$ is not tight. The next result
demonstrates the tightness of the order of magnitude of the
bound $\sqrt{2k\log d}$ for large values of
$d$.\footnote{I wish to thank Benjamin Weiss for raising the
question of the tightness of the factor $\sqrt{\log d}$ in the bound,
and demonstrating for each positive $\ell$ the existence of a simple martingale of
probabilities $p_0^{\ell}$ on a set with $2^{\ell}$ elements and
with variation $\ell$. Specifically, starting with the uniform probability, in
each stage half of the non-zero probabilities (each half equally
likely) move to zero, and the other half double their probabilities.
Therefore, for each fixed $\alpha>0$ there is a positive constant
$0<C(\alpha)$ ($\to_{\alpha\to 0+}0$) such that for $k$ and $d$ with
$\alpha \leq \frac{\log_2 d}{k} \leq 1$, $V(k,d)\geq [\log_2 d]\geq
C(\alpha)\sqrt{k\log d}$.} We have

\begin{theorem}\label{thm:2} There is a positive constant $C>0$ such that for every
$k$
 and $d$ with $d\leq 2^k$ there is $p_0^k\in \mcM_k(X,p_0)$ with
$| X|=d$ such that
\[V(p_0^k)\geq C\sqrt{k\log d }\,.\]
\end{theorem}

Bounds of the variation of martingales of probabilities are useful
in the study of repeated games with incomplete
information  \cite{aummas95}. In a two-person
zero-sum repeated game with incomplete information on one side
(henceforth, RGII-OS) the players play repeatedly the same stage game
$G$. However, the game depends on a state $x\in X$ known only to
player 1 (P1) and $x$ is chosen according to a probability $p\in
\Delta(X)$ that is commonly known. In the course of the game player
2 (P2) may learn information about $x$ only from past actions of
player 1.

Formally, a RGII-OS $\Gamma$ is defined by a state space $X$, a
probability $p\in \Delta(X)$, finite sets of stage actions, $I$ for
P1 and $J$ for P2,  and for every $x\in X$ we have a two-person zero-sum
$I\times J$ matrix game $G^x$. We write $\Gamma=\langle
X,p,I,J,G\rangle$, where $G$ stands for the list of
matrix games $(G^x)_{x\in X}$. The $(i,j)$-th entry of $G^x$,
denoted $G^x_{i,j}$, is the payoff from P2 to P1 when in state $x$
the players play the action pair $(i,j)$.

The $k$-stage repeated game, denoted $\Gamma_k(p)$, or $\Gamma_k$
for short, is played as follows. Nature chooses $x\in X$ according
to the probability $p$. P1 is informed of nature's choice $x$, but P2 is
not. At stage $1\leq t\leq k$, P1 chooses $i_t\in I$ and
simultaneously P2 chooses $j_t\in J$ (and these choices are observed
by the players following the play in stage $t$). The choice of $i_t$
may depend on $x,i_1,j_1,\ldots,i_{t-1},j_{t-1}$ (which is the
information of P1 before the play at stage $t$) and the choice of
$j_t$ may depend on $i_1,j_1,\ldots,i_{t-1},j_{t-1}$ (which is the
information of P2 before the play at stage $t$).

A  pair of strategies $\s$ of P1 and $\tau$ of P2 (together with the
initial probability $p$) define a probability distribution
$P^p_{\s,\tau}$, or $P_{\s,\tau}$ for short, on the space of plays
$x,i_1,j_1,\ldots,i_k,j_k$, and thus on the stream of payoffs
$g_t:=G^x_{i_t,j_t}$. The (normalized) payoff of the $k$-stage
repeated game is the average of the payoffs in the $k$-stages of the
game, namely, $\bar{g}_k=\frac1k\sum_{t=1}^kg_t$. The minmax value
of $\Gamma_k(p)$ is $v_k(p):=\max_{\s}\min_{\tau}
E_{\s,\tau}\bar{g}_k$, where $E_{\s,\tau}$ stands for the expectation
with respect to the probability $P^p_{\s,\tau}$, the maximum is over all mixed
(or behavioral) strategies $\s$ of P1, and the minimum is over all mixed
(or behavioral) strategies $\tau$  of P2.

For fixed components $\langle X,I,J,G\rangle$, the minmax value of the matrix game $\sum_xq(x)G^x$
 is a function of $q\in \Delta (X)$ and is denoted $u(q)$.
 The least concave function on $\Delta(X)$ that is greater than or equal to $u$ (``smallest concave majorant'') is
denoted $\cav\,u$. Aumann and Maschler \cite{aummas95} proved that
$v_k(p)\geq (\cav\, u)\,(p)$ and that $v_k(p)$ converges to
$(\cav\,u)\,(p)$ as $k\to \infty$ . Moreover, \cite{aummas95} shows
that the bound of the variation of the martingale of probabilities
bounds the (nonnegative) difference $v_k(p)- (\cav\, u)\,(p)$.
Explicitly, if $\|G\|:=\max_{x,i,j}|G^x_{i,j}|$, we have
\beqn\label{eqn:valuebound} v_k(p)- (\cav\, u)\,(p)\leq \|G\|
V(k,p)/k. \eeqn
Inequality (\ref{eqn:valuebound}) yields on the one hand
a rate of convergence of $v_k(p)$, and on the other hand enables us
to approximate the value $v_k(p)$ for a specific $k$ and a specific
game. The classical bound of $V(k,p)$ that is used in
\cite{aummas95, merzam71} and in subsequent works is
\[V(k,p)\leq V(k,d)\leq \sqrt{k(d-1)}\,.\]
For $d>k$ this bound is not useful. Theorem \ref{thm:martingalevariationbound} provides an effective
bound when $d$ is subexponential in $k$, namely, when $\log d=o(k)$,
or, more generally, when $H(p)/k$ is small. Applying the bound in
Theorem \ref{thm:martingalevariationbound} to the inequality
(\ref{eqn:valuebound}) implies that \beqn\label{eqn:valuebound1}
v_k(p)- (\cav\, u)\,(p)\leq \|G\| \sqrt{\frac{2\log d}{k}}\,. \eeqn

One may wonder if the order of magnitude of the bound in
(\ref{eqn:valuebound1}) is tight. We have
\begin{theorem}\label{thm:rgiiierrorboundfrom below}
There is a positive constant $C$ such that for every $k$ and $d$
with $d\leq 2^k$ there is a repeated game $\Gamma$ with incomplete
information on one side  with ($\|G\|>0$ and) $d$ states
such that \beqn v_k(p)- (\cav\, u)\,(p)\geq C\|G\| V(k,d)/k.\eeqn
\end{theorem}

\section{Proofs}

{\em Proof of Theorem \ref{thm:martingalevariationbound}.} Let $p_0^k$ be $(\mcH_t)_t$-adapted; that is,
$p_t$ is measurable with respect to the $\s$-algebra $\mcH_t\subset
\mcH_{t+1}$. Without loss of generality we can assume that $\mcH_t$ are finite, namely, algebras. (Indeed, if $p_t$ is measurable
with respect to the $\s$-algebra $\mcH_t\subset \mcH_{t+1}$, one replaces
$\mcH_t$ with an algebra $\mcH^*_t\subset \mcH_{t+1}^*$ such that
$\|E(p_t\mid \mcH_t^*)-p_t\|\leq \ve/k$, and replaces $p_t$ with
$\hat{p}_t:=E(p_t\mid \mcH_t^*)$ ($=E(p_k\mid \mcH_t^*)$). Note that
$\sum_{t=1}^k \|\hat{p}_{t}-\hat{p}_{t-1}\|+2\ve\geq \sum_{t=1}^k
\|p_t-p_{t-1}\|$.) In
that case we can assume that: (1) $P$ is a probability on the product
$X\times (\times_{t=0}^k A_t)$, where $A_t$ are finite sets (e.g.,
the atoms of the algebra $\mcH_t$); (2) $(x,a_0,a_1,\ldots,a_k)$ is
a vector of random variables having distribution $P$; and (3) $p_t$
is the conditional distribution of $x$ given $a_0,\ldots,a_t$. Let
$P_t$ be the conditional (joint) distribution of $(x,a_t)$ given
$(a_0,\ldots,a_{t-1})$, $P_{tX}$ its marginal on $X$, and $P_{tA_t}$
its marginal on $A_t$. Let $P_{tX}\otimes P_{tA_t}$ denote the product distribution on $X\times A_t$, i.e.,
$P_{tX}\otimes P_{tA_t}(x,a_t)=P_{tX}(x)P_{tA_t}(a_t)$.
By Pinsker's inequality (see, e.g., \cite[p.
300]{covtho}), we have
\[\|P_t-P_{tX}\otimes P_{tA_t}\|\leq \sqrt{2}\sqrt{D(P_t\|P_{tX}\otimes P_{tA_t})}\,,\]
where for two probabilities $P$ and $Q$ on a finite (or
countable) set $Y$, $\|P-Q\|=\sum_y|P(y)-Q(y)|$ and $D(P\|Q)=\sum_{y\in Y}P(y)\log\frac{P(y)}{Q(y)}$
(where $\log$ denotes the natural logarithm and $0\log 0=0$).

Let $H_{P_t}(x):=-\sum_x P_t(x)\log P_t(x)$, $H_{P_t}(a_t)$, and
$H_{P_t}(x, a_t)$, denote the entropy of the random variables $x$,
$a_t$, and $(x,a_t)$, where $(x,a_t)$ has distribution $P_t$, and
$H_{P_t}(x\mid a_t):=H_{P_t}(x, a_t)-H_{P_t}(a_t)$. A straightforward computation yields 
$D(P_t\|P_{tX}\otimes
P_{tA_t})=H_{P_t}(x)-H_{P_t}(x\mid a_t)$.
Therefore,
\beqn \|P_t-P_{tX}\otimes P_{tA_t}\|\leq \sqrt{2}
\sqrt{H_{P_t}(x)-H_{P_t}(x\mid a_t)}\,. \eeqn
Note that $P_t$ is a random variable, which is a function of
$a_0,\ldots,a_{t-1}$, and therefore, by the properties of conditional
entropy, $E_PH_{P_t}(x)=H_P(x\mid a_0,\ldots,a_{t-1})$ (where $E_P$
denotes the expectation with respect to the probability distribution $P$) and
$E_PH_{P_t}(x\mid a_t)=H_P(x\mid a_0,\ldots,a_{t-1},a_t)$. Therefore,
\[E_P \left(H_{P_t}(x)-H_{P_t}(x\mid a_t)\right)=H_P(x\mid
a_0,\ldots,a_{t-1})-H_P(x\mid a_0,\ldots,a_{t-1},a_t).\] As the
square root is a concave function we have, by Jensen's inequality,
\[E_P \|P_t-P_{tX}\otimes P_{tA_t}\|\leq \sqrt{2}\sqrt{H_P(x\mid
a_0,\ldots,a_{t-1})-H_P(x\mid a_0,\ldots,a_{t-1},a_t)}.\]

As %the conditional expectation
$E_P(\|p_t-p_{t-1}\| \mid \mcH_{t-1})$ equals $\sum_{a\in
A_t}P_{tA_t}(a)\sum_x|\frac{P_t(x,a)}{P_{tA_t}(a)}-P_{tX}(x)|=\sum_{a\in
A_t}\sum_x|{P_t(x,a)}-P_{tA_t}(a)P_{tX}(x)|=\|P_t-P_{tX}\otimes
P_{tA_t}\|$, we deduce that $E_P\|p_t-p_{t-1}\|=
E_P\|P_t-P_{tX}\otimes P_{tA_t}\|$ and therefore by substituting
$E_P\|p_t-p_{t-1}\|$ for $E_P\|P_t-P_{tX}\otimes P_{tA_t}\|$ we get
\[E_P \|p_t-p_{t-1}\|\leq \sqrt{2}\sqrt{H_P(x\mid
a_0,\ldots,a_{t-1})-H_P(x\mid a_0,\ldots,a_{t-1},a_t)}\,.\]

As the square root is a concave function,  using Jensen's inequality
and the equality and inequality  $\sum_{t=1}^k\left(H(x\mid
a_0,\ldots,a_{t-1})-H(x\mid a_0,\ldots,a_t)\right)=H(x)-H(x\mid
a_0,\ldots,a_k)\leq H(x)$, we have
\[E \sum_{t=1}^k \|p_{t}-p_{t-1}\| \leq \sqrt{2k}\sqrt{H(x)}\leq \sqrt{2k}\sqrt{\log d}\,.\]
This completes the proof of Theorem
\ref{thm:martingalevariationbound}.
\endproof

{\em Proof of Theorem \ref{thm:2}.} Note that $V(k,d)$ is
monotonic increasing in $d$ and $k$, and there is a positive
constant $C_1>0$ such that $V(k,2)\geq C_1\sqrt{k}$.

If $p_0^{k_1}$ and $q_0^{k_2}$ are two martingales with total
variation $V_1$ and $V_2$, respectively, then $p_0\otimes
q_0,\ldots,p_{k_1}\otimes q_0$ is a martingale with total variation
$V_1$ and $p_{k_1}\otimes q_0,p_{k_1}\otimes q_1,\ldots,p_{k_1}\otimes
q_{k_2}$ is a martingale with total variation $V_2$ and therefore
$p_0\otimes q_0,\ldots,p_{k_1}\otimes q_0,p_{k_1}\otimes
q_1,\ldots,p_{k_1}\otimes q_{k_2}$ is a martingale with total variation
$V_1+V_2$. Therefore, %This completes the proof of (\ref{eqn:6}).
%A simple inequality is
%
\beqn\label{eqn:6} V(k_1,p)+V(k_2,q)\leq V(k_1+k_2,p\otimes q), \eeqn
from which it follows that
\beqn\label{eqn:addedvariation} V(k_1,d_1)+V(k_2,d_2)\leq
V(k_1+k_2,d_1d_2).
\eeqn

Inequality (\ref{eqn:addedvariation}) implies that if $k$ is a
multiple of $\ell$ we have $V(k,2^{\ell})\geq \ell V(k/\ell,2)\geq
\ell C_1\sqrt{k/\ell}=C_1\sqrt{k\ell}$. Note that for every $k$ and
$2\leq d\leq 2^k$ there is $k\geq k_1>k/2$ that is a multiple of
$\ell=[\log_2 d]\geq (\log_2 d)/2$ (where $[x]$ is the largest
integer $\leq x$), and therefore $V(k,d)\geq V(k_1,2^{\ell})\geq
C_1\sqrt{k_1\ell}\geq  C_1/2\sqrt{k\log_2 d}$. This completes the
proof of Theorem \ref{thm:2}.
\endproof

{\em Proof of Theorem \ref{thm:rgiiierrorboundfrom below}.}
Given two repeated games with incomplete information on one side,
$\Gamma^1=\langle X_1,p_1,I_1,J_1,G^1\rangle$ and $\Gamma^2=\langle
X_2,p_2,I_2,J_2,G^2\rangle$, we define the game
$\Gamma=\Gamma_1\otimes \Gamma_2$ by
\[\Gamma =\langle X=X_1\times
X_2,p=p_1\otimes p_2,I=I_1\times I_2\times \{1,2\},J=J_1\times
J_2,G\rangle\,,\]
where for $x=(x_1,x_2)$, $i=(i^1,i^2,b)$, and
$j=(j^1,j^2)\in J$,
\[G^x_{i,j}=G^{x_b}_{i^b,j^b}\,,\]
where $G^{x_b}$ stands for the more explicit $G^{b,x_b}$.
Note that $\|G\|=\max (\|G^1\|,\|G^2\|)$.

A possible helpful interpretation of $\Gamma$ is that nature chooses
a pair $x_1\in X_1$ and $x_2\in X_2$, equivalently a pair of games
$G^{x_1}$ and $G^{x_2}$,  according to the product probability
$p_1\otimes p_2$. P1 is informed of the choice $(G^{x_1},G^{x_2})$
of nature, but P2 is not. In each stage of the repeated game, both
players select strategies for the first and for the second game, and
P1 chooses in addition which one of the two games determines the
stage payoff.

 As a function of
$i=(i^1,i^2,b)$, for each fixed $b=1,2$, the payoff function
$G^x_{i,j}$ does not depend on the coordinate $i^c$ for $c\neq b$.
Therefore we can replace the set $I$ (which has $2|I_1||I_2|$
elements) of stage actions of P1 in the repeated game $\Gamma$ with
the disjoint union of $I_1$ and $I_2$.

Note that if $v^b_k$ and $v_k$ stand for the (normalized) values of
the $k$-stage repeated games $\Gamma^b$ and $\Gamma$,  then
\[v_{k_1+k_2}\geq \frac{k_1v^1_{k_1}+k_2 v^2_{k_2}}{k_1+k_2}.\]
Indeed, P1 can play $b_t=1$ in stages $t=1,\ldots,k_1$ and $b_t=2$
in stages $t=k_1+1,\ldots,k_1+k_2$, and the first coordinates
$i^1_t$ of $i_t$ follow, in stages $t=1,\ldots,k_1$, an optimal
strategy of P1 in $\Gamma^1_{k_1}(p_1)$, and the second coordinates
$i^2_t$ of $i_t$ follow, in stages $t=k_1+1,\ldots,k_1+k_2$, an
optimal strategy of P1 in $\Gamma^2_{k_2}(p_2)$.

For $\ell>2$ and a sequence $\Gamma^1=\langle
X_1,p_1,I_1,J_1,G^1\rangle,\ldots,\Gamma^{\ell}=\langle
X_\ell,p_\ell,I_\ell,J_\ell,G^\ell\rangle$ of RGII-OS,  we define by
induction on $\ell$ the game $\Gamma=\otimes_{b=1}^\ell\Gamma^b$ by
 $\Gamma
=\left(\otimes_{b=1}^{\ell-1}\Gamma^b\right)\otimes \Gamma^\ell$.

If $v^b_k$, respectively $v_k$, denotes the normalized value of the
$k$-stage repeated game $\Gamma^b_k(p_b)$, respectively
$\Gamma_k(\otimes_{b=1}^{\ell}p_b)$, and $k=k_1+\ldots+k_\ell$, then
\[v_k\geq \frac{\sum_{b=1}^\ell k_bv^b_{k_b}}{k}.\]
Note that a stage action of P1 in $\Gamma$ is a list of stage
actions $i^1,\ldots,i^\ell$ (with $i^b\in I_b$) and a number $b$ (with $1\leq
b\leq \ell$). However, given $b$, the payoff depends only on the
coordinate $i^b$ of the stage actions. Therefore we can replace the
stage actions of P1 in $\Gamma$ with the disjoint union of the
action sets $I_b$, and so with a set of size $\sum_b |I_b|$.

Consider the example of the RGII-OS $\Gamma^z=\langle
X=\{0,1\},(1/2,1/2), I,J,G\rangle$, introduced by Zamir
\cite[Section 3]{zam72}. The set of states is $X=\{0,1\}$, and
players' action sets are $I=\{0,1\}$ for P1, and $J=\{0,1\}$ for P2.
The two payoff matrices are $G^0$ and $G^1$: $G^0_{0,0}=3$,
$G^0_{0,1}=-1$, $G^1_{0,0}=2=-G^1_{0,1}$, and
$G^*_{i,j}=-G^*_{1-i,j}$. Let $v^z_k$ denote the normalized value of
the $k$-stage repeated game $\Gamma^z$. Zamir \cite{zam72} shows
that $\lim_n v^z_n=0$ and $v^z_k\geq C_1/\sqrt{k}$, where $C_1>0$ is
a positive constant.
%
\begin{comment}In addition, there is a constant
$C_2$ and a strategy of P2 in the infinitely repeated that
guarantees that for every strategy $\s$ of P1 the expected average
payoff in the first $n$ stages is $\leq C_2/\sqrt{n}$.
\end{comment}

Consider the RGII-OS $\Gamma=\otimes_{b=1}^{\ell}\Gamma^z$, and let
$v_k$ denote the normalized value of the $k$-stage repeated game
$\Gamma$. It follows that
\[v_k\geq \max\left\{\frac{\sum_{b=1}^\ell k_bv^z_{k_b}}{k}:
k_b\geq 0 \mbox{ and }\sum_{b=1}^\ell k_b=k\right\},\] and therefore
if $k$ is a multiple of $\ell$ we can take $k_b=k/\ell$ and
therefore
\[v_k\geq v_{k/\ell}^z\geq C_1\sqrt{\ell/k}\,.\]
For an arbitrary $k$ and $d\leq 2^k$, there is $k\geq k_1>k/2$ that
is a multiple of $\ell:=[\log_2 d]$ ($\geq (\log_2 d)/2$). As P1 can
play $(1/2,1/2)$ in the last $k-k_1$ stages, $kv_k\geq k_1v_{k_1}$,
and thus $kv_k\geq  C_1\sqrt{k_1}\sqrt{\ell}$.
Therefore, $v_k\geq \frac{C_1}{2\sqrt{k}}\sqrt{\log d}$.

Finally, the existence of an optimal strategy of P2 in the
infinitely repeated game $\Gamma^z$ (or a direct computation of the
function $u(p)$, the minmax value of the game $\sum_x p(x)G^x$, for
the game $\Gamma$) yields $\lim_n v_n=0$. Note that the stage
payoffs of the RGII-OS $\Gamma$ are bounded by $3$ (independent of
the number of factors $\ell$). Altogether, we have constructed for
each $k$ and $d\leq 2^k$ a repeated game $\Gamma=\langle
X,p,I,J,G\rangle$ with $|X|\leq d$, equivalently $|X|=d$ ($|I|\leq
2\log d$) and $\|G\|=3$, and
\[v_k-\lim_{n\to \infty}v_n \geq C_1/2\sqrt{\log d}/\sqrt{k}\,.\] This
completes the proof of Theorem \ref{thm:rgiiierrorboundfrom
below}.
\endproof

\section{Remarks}
\subsection{Comments on the proof of Theorem 1.}
Our proof of Theorem 1 relies on Pinsker's inequality, and it uses information-theoretic tools. In fact, the information-theoretic intuition has led us to the result and its proof.
However, readers unfamiliar with the information-theoretic concepts may find the proof obscure.
The following is an alternative derivation (which disguises the use of the information-theoretic techniques) and uses classical martingale theory techniques.
First, note that Pinsker's inequality implies that if $Z$ and $Y$ are two nonnegative random variables with $EZ=EY$, then $E|Z-Y|\leq \sqrt{2EZ}\sqrt{EZ\log Z-EZ\log Y}$.
In particular, if $E(Z\mid Y)=Y$ (e.g., when $Y$ is the constant random variable $Y=EZ$), then $E Z\log Y=E(E(Z\log Y \mid Y))=EY\log Y$, and therefore
\begin{equation}\label{eqn:pinskervariant}E|Z-Y|\leq \sqrt{2EZ}\sqrt{EZ\log Z-EZ\log Y}=\sqrt{2EZ}\sqrt{EZ\log Z-EY\log Y}.\end{equation}
Inequality (\ref{eqn:pinskervariant}), which is equivalent to Pinsker's inequality, can obviously be proved directly.\footnote{ I wish to thank Stanislaw Kwapien for suggesting a proof that avoids the information-theoretic techniques, and communicating a simple analytical proof of the above displayed version of Pinsker's inequality.} The continuation of the proof avoids the (explicit) use of information-theoretic techniques.

It follows from (\ref{eqn:pinskervariant}) that
if $Z_0,\dots, Z_k$ is a martingale of nonnegative random variables, then \[\sum_{t=1}^k |Z_t-Z_{t-1}|\leq \sqrt{2EZ_0}\sum_{t=1}^k \sqrt{EZ_tlog Z_t-EZ_{t-1}\log Z_{t-1}},\] which
by Jensen's inequality, the concavity of the square root, and the telescopic feature of the series  $EZ_t\log Z_t-EZ_{t-1}\log Z_{t-1}$, is \[\leq \sqrt{2kEZ_0}\sqrt{EZ_klog Z_k-EZ_{0}\log Z_{0}}.\] Therefore, if $p_0,\ldots,p_k$ is a martingale with values in $\mathbb{R}^X_+$
we have
\[\sum_{t=1}^k E\|p_t-p_{t-1}\|_1\leq \sum_{x\in X}\sqrt{2kE p_0(x)}\,\sqrt{Ep_k(x)\log p_k(x)- Ep_0(x)\log p_0(x)}.\]
By the Schwartz inequality we obtain that
\[\sum_{t=1}^k E\|p_t-p_{t-1}\|\leq \sqrt{2kE \sum_{x\in X}p_0(x)}\,\sqrt{E\sum_{x\in X}\left(p_k(x)\log p_k(x)- p_0(x)\log p_0(x)\right)}\,.\]
If $p_k(x)\leq M(x)$, then by the convexity of $q\log q$ we have $Ep_k(x)\log p_k(x)\leq Ep_0(x)\log M(x)$, and then
\[\sum_{t=1}^k E\|p_t-p_{t-1}\|_1\leq \sqrt{2kE \sum_{x\in X}p_0(x)}\sqrt{\sum_{x\in X}- Ep_0(x)\log (p_0(x)/ M(x))}\,.\]
If $p_k(x)\leq 1$, then $p_k(x)\log p_k(x)\leq 0$, and therefore
\[\sum_{t=1}^k E\|p_t-p_{t-1}\|_1\leq \sqrt{2kE \sum_{x\in X}p_0(x)}\sqrt{\sum_{x\in X}- Ep_0(x)\log p_0(x)}\,.\]
We conclude that
if $\sum_x p_0(x)=1$, then $\sum_{t=1}^k E\|p_t-p_{t-1}\|_1\leq \sqrt{2k}\sqrt{EH(p_0)}$.

\subsection{The variation of a martingale of probabilities over a countable set.}
It is of interest to find a necessary and sufficient
condition for a distribution $p$ on a countable set $X$ for
$\sup_k\frac{1}{\sqrt{k}}V(k,p)<\infty$. We remark here on the
sufficient conditions derived from the classical method and our
method of bounding the variation of martingales of probabilities.

The classical bound of the variation of a martingale $p_0^k$ is
obtained by bounding, for each fixed $x\in X$, the expectation
variation $E\|y(x)\|_1$, where $y(x)\in \mathbb{R}^k$ is the vector
of martingale differences $(p_1(x)-p_0(x),\ldots,p_k(x)-p_{k-1}(x))$
(thus $\|y(x)\|_1=\sum_{t=1}^k |p_t(x)-p_{t-1}(x)|$), and summing
over all $x\in X$. Assuming without loss of generality that $p_0$ is a constant $p\in
\Delta(X)$, we have (by the Cauchy--Schwartz inequality)
$\|y(x)\|_1\leq \sqrt{k}\|y(x)\|_2$, and, therefore, by Jensen's
inequality, $E \|y(x)\|_1\leq \sqrt{k}\sqrt{E\|y(x)\|_2^2}$, which
by the martingale property is $\leq
\sqrt{k}\sqrt{E((p_k(x))^2-(p_0(x))^2)}\leq
\sqrt{k}\sqrt{E((p_k(x))-(p_0(x))^2)}=\sqrt{k}\sqrt{p_0(x)-(p_0(x))^2}$.
Therefore, if $p\in \Delta(X)$ and $X$ is countable, the classical
method yields that $\sup_k\frac{1}{\sqrt{k}}V(k,p)<\infty$ whenever
$\sum_x \sqrt{p(x)}<\infty$. As $-q\log q =o(\sqrt{q})$ as $q\to
0+$, the condition  $\sum_x \sqrt{p(x)}<\infty$ implies that
$H(p)=-\sum_x p(x)\log p(x)<\infty$. Obviously, there are
probabilities $p$ over a countable set $X$ such that $H(p)<\infty$
but $\sum_x \sqrt{p(x)}=\infty$. Therefore our bound provides a
strictly sharper sufficient condition, $H(p)<\infty$, for
$\sup_k\frac{1}{\sqrt{k}}V(k,p)<\infty$, compared to the one derived
by using the classical method.

\subsection{The asymptotic behavior of $V(k,d)$.}
The  asymptotic behavior of $V(k,d)$ deserves further
study. \cite{merzam77} proves that ${V(k,2)}/{\sqrt{k}}$ converges
as $k\to \infty$ to $\sqrt{\frac{2}{\pi}}$. It is of interest to
find a corresponding limit theorem for ${V(k,d)}/{\sqrt{k\log d}}$
as $2^k\geq d\to \infty$. The above-mentioned result of
\cite{merzam77} together with our construction in the proof of
Theorem 2 yields that the $\liminf$ of ${V(k,d)}/{\sqrt{k\log_2 d}}$
is $\geq \sqrt{\frac{2}{\pi}}$ as $\frac{\log d}{k}+ 1/d\to 0$
(namely, as $\log d =o(k)$ and $d\to \infty$).

\begin{comment}
The entropy $H(p)$ of the distribution $p$ is an insufficient
parameter for the study of the asymptotic behavior of $V(k,p)$ ,
even for large values of $H(p)$; namely, there are sequences of
probabilities over finite sets, $p^k$ and $q^k$, with
$H(p^k)=H(q^k)\to_{k\to \infty} \infty$ for which
$\frac{V(k,p^k)-V(k,q^k)}{\sqrt{k H(p^k)}}$ is  bounded away from 0
for sufficiently large $k$. Moreover, for every sequence $H_k>0$
with $k\geq H_k\to_{k\to \infty}\infty$, one can find such sequences
$p^k$ and $q^k$ with $H_k=H(q^k)$. For example, for each $k\geq 2$
let $\alpha_k$ be the solution of
$\frac{H(\alpha_k,1-\alpha_k)}{\alpha_k}=H_k$ in the interval
$0<\alpha_k<1$. Note that $\alpha_k\to 0$ as $k\to \infty$. Define
the probability $q^k$ on the nonnegative integers by
$q^k(i)=\alpha_k(1-\alpha_k)^i$. Note that the entropy of $q^k$,
$H(q^k)$, equals $H_k$. By Theorem \ref{thm:}, there is a positive
constant $C_1>0$ such that for every $k$ there is a probability
$p^k$ on (a finite set) such that $H(p^k)=H_k$ and $V(k,p^k)\geq
C_1\sqrt{kH_k}$. Note that $V(k,q^k)\leq
\sum_{i=1}^{\infty}V(k,(q^k(i),1-q^k(i)))\leq
\sum_{i=1}^{\infty}\sqrt{kq^k(i)}\leq \sqrt{k/\alpha_k}$. It follows
that
\[\frac{V(k,p^k)-V(k,q^k)}{\sqrt{k H(p^k)}}\geq C_1-\sqrt{\frac{1}{\alpha_k H_k}}.\]
\end{comment}

\subsection{Repeated games with incomplete information.}
The proof of Theorem \ref{thm:rgiiierrorboundfrom below} constructs for each $d\leq 2^k$ a RGII-OS $\Gamma=\langle
X,p,I,J,G\rangle$ with $|X|=d$ and $v_k\geq \lim_nv_n+
C\sqrt{k^{-1}\log d}$, where $0<C=O(\|G\|)$, and, in addition, $|I|=
O(\log d)$ and $|J|=O(d)$. We have not tried to minimize the order
of magnitude of the number of elements of $I$ and $J$. It is however
impossible to construct such an example with bounded $|I|$ and
$|J|$. Indeed, in a forthcoming note we will show that for every RGII-OS
$\Gamma=\langle X,p,I,J,G\rangle$ we have $v_k\leq \lim_n
v_n+\|G\|_* \sqrt{\frac{2|I\times J|\log k}{k}}$, where $\|G\|_*:=2E_p(\max_{i,j}G^x_{i,j}-\min_{i,j}G^x_{i,j})$. Therefore the
inequality $v_k\geq \lim_nv_n+ C\sqrt{k^{-1}\log d}$ is possible
only if $|I\times J|\geq C\frac{\log d}{2\log k}$.

\end{document}